\documentclass[a4paper,11pt]{article}
\input epsf
\usepackage{latexsym}
\usepackage{amssymb}
\newtheorem{thm}{Th\'eor\`eme}[section]
\newtheorem{lem}[thm]{Lemma}

\newtheorem{prop}[thm]{Proposition}

\newtheorem{rem}[thm]{Remarque}

\input epsf

\title{La factorisation de  $x^{p^l}-x\in {\mathbb F}_p[x]$ 
selon la trace}
\author{Roland Bacher}
\date{}

\begin{document}
\maketitle


{\bf R\'esum\'e:} Nous d\'ecrivons quelques factorisations de
polyn\^omes sur des corps finis. Ces factorisation sont
li\'ees \`a la trace, aux compositions de polyn\^omes et
aux coefficients binomiaux. Comme cons\'equence nous obtenons
la description des polyn\^omes irr\'eductibles $Q\in{\mathbb F}_2[x]$
tels que les polyn\^omes $Q(1+x+x^2)$ (ou $Q(x+x^2)$) sont \'egalement
irr\'eductibles. 

{\bf Abstract:} We present a few factorizations of polynomials over
finite fields. These factorizations are related to traces,
compositions of polynomials and binomial coefficients. As a corollary
we obtain a description of all irreducible polynomials  $Q\in
{\mathbb F}_2[x]$ such that $Q(1+x+x^2)$ (or $Q(x+x^2)$) 
remain irreducible.

\section{R\'esultats principaux}

Pour $p$ un nombre premier, les racines du polyn\^ome $x^{p^l}-x\in
{\mathbb F}_p[x]$ sont les \'el\'ements du corps
fini ${\mathbb F}_{p^l}={\mathbb F}_p[x]/(R)$ \`a $q=p^l$ 
\'el\'ements o\`u $R\in{\mathbb F}_p[x]$ est un
polyn\^ome irr\'eductible de degr\'e $l$. Rappelons la d\'efinition
de la {\it trace
relative} $\hbox{tr}_l(\rho)=\sum_{k=0}^{l-1}\rho^{p^k}\in {\mathbb F}_p$
d'un \'el\'ement $\rho\in{\mathbb F}_{p^l}$. C'est donc la trace 
de l'application $x\longmapsto \rho x$, consid\'er\'ee comme
endomorphisme ${\mathbb F}_p-$lin\'eaire de ${\mathbb F}_{p^l}$.
Pour $F=\sum_{j=0}^d\alpha_jx^j\in{\mathbb F}_p[x]$ 
un polyn\^ome irr\'eductible de degr\'e $d$ divisant $l$ nous posons
$\hbox{tr}_l(F)=\hbox{tr}_l(\rho)=-\frac{l}{d}\frac{\alpha_{d-1}}
{\alpha_d}$ pour d\'efinir sa trace relative o\`u
$\rho\in {\mathbb F}_{p^l}$ est une racine de $F$.

Dans la suite, on identifiera g\'en\'eralement deux diviseurs 
$F,G$ d'un polyn\^ome $P\in{\mathbb F}_p[x]$ si $G=\lambda F$ pour 
$\lambda\in {\mathbb F}_p^*$. L'ensemble des diviseurs d'un
polyn\^ome $P$ sera d\'efini comme l'ensemble
des polyn\^omes moniques divisant $P$.

Le r\'esultat suivant est, au moins partiellement, bien connu, 
voir par exemple
la proposition 3.4.7 dans \cite{C}.

\begin{thm} \label{tracefact} (i) Le polyn\^ome 
$$-\alpha+\sum_{k=0}^{l-1}x^{p^k}\in{\mathbb F}_p[x]$$
est le produit de tous les facteurs irr\'eductibles $F$ divisant
$x^{p^l}-x$ dont la trace relative est $\hbox{tr}_l(F)=\alpha$.

\ \ (ii) Pour $\alpha\in{\mathbb F}_p$ et $d$ un diviseur de $l=df$,
  le polyn\^ome 
$$-\alpha+\sum_{k=0}^{f-1}x^{p^{dk}}\in{\mathbb F}_p[x]$$
divise $-d\alpha+\sum_{k=0}^{l-1}x^{p^k}\in {\mathbb F}_p[x]$.
\end{thm}

Notre preuve d\'emontre une l\'eg\`ere g\'en\'eralisation du 
th\'eor\`eme \ref{tracefact} : On peut rempla\c cer le corps 
primaire ${\mathbb F}_p$ par le corps fini ${\mathbb F}_q$
\`a $q=p^e$ \'el\'ements.

La factorisation partielle
$$X^{p^l}-X=\prod_{\alpha=0}^{p-1}\left(-\alpha+\sum_{k=0}^{l-1}
X^{p^k}\right)\in {\mathbb F}_p[X]$$
est l'ingr\'edient principal de l'algorithme 3.4.8 dans \cite{C}
permettant la factorisation dans
${\mathbb F}_2[x]$. La caract\'erisation des diviseurs irr\'eductibles
de $-\alpha+\sum_{k=0}^{l-1}x^{p^k}\in{\mathbb F}_p[x]$ en fonction de
leurs traces simplifie l\'eg\`erement l'\'etude de cet algorithme. 

Avant de passer en caract\'eristique $2$, mentionnons encore le
r\'esultat suivant (probablement bien connu) qui est valable pour
un corps commutatif quelconque.

\begin{prop} \label{decomp} Consid\'erons deux 
polyn\^omes $Q,R\in K[x]$ \`a coefficients dans un
corps commutatif $K$. Supposons $Q$ irr\'eductible.
Alors son degr\'e $\hbox{deg}(Q)$ divise le degr\'e 
de tout facteur irr\'eductible 
$F\in K[x]$ du polyn\^ome compos\'e $Q\circ R\in K[x]$.
\end{prop}

La proposition \ref{decomp} associe donc \`a une paire de 
polyn\^omes $Q,R\in K[x]$ avec 
$Q$ irr\'eductible 
une partition $\sum_{j=1}^a r_j$ du degr\'e de $R$ 
obtenue en consid\'erant les degr\'es 
$(r_j\ \hbox{deg}(Q))=\hbox{deg}(F_j)$ des $a$
facteurs irr\'eductibles (pas n\'ecessairement distincts)
du polyn\^ome compos\'e $Q\circ R=F_1\cdots F_a$. Un cas particulier
amusant est la partition associ\'ee \`a $Q\circ Q$ pour $Q\in 
K[x]$ irr\'eductible. Y-a-t-il des restrictions sur les partitions
obtenues \`a partir de $Q\circ Q$? (Cela semble \^etre le
cas sur le corps ${\mathbb F}_2$: 
Les partitions $l=1+1+\dots+1$ n'apparaissent pas pour 
$Q\in{\mathbb F}_2[x]$ irr\'eductible de degr\'e $l\in\{2,\dots,8\}$.)
On peut \'egalement se poser des questions sur les fr\'equences
(asymptotiques) des partitions associ\'ees \`a $Q\circ Q$
(respectivement $Q\circ R$)
pour $Q$ irr\'eductible et $Q$ (respectivement $R$) de degr\'e grand etc. 

Pour $P,Q\in{\mathbb F}_q[x]$ \`a coefficients dans le corps fini
\`a $q=p^e$ \'el\'ements, $Q$ irr\'eductible de degr\'e $l$,
le nombre de facteurs irr\'eductibles distincts $F\in{\mathbb F}_q[x]$ 
divisant $Q\circ R$ de degr\'e un diviseur de $dl$ est \'evidemment
donn\'e par le nombre de facteurs irr\'eductible du plus grand 
diviseur commun entre $Q\circ R$ et $x^{q^{dl}}-x\pmod Q\circ R$.

Le premier cas non-trivial, obtenu en choisissant $R$ parmi les
polyn\^omes
$x+x^2,1+x+x^2\in{\mathbb F}_2$ de degr\'e $2$ \`a coefficients
dans ${\mathbb F}_2$ admet une description plus simple,
donn\'ee par le r\'esultat suivant. (Les cas 
$Q(x^2)=(Q(x))^2$ et $Q(1+x^2)=(Q(1+x))^2\in{\mathbb F}_2[x]$ sont sans
int\'er\^et.)

\begin{thm} \label{decomp2} 
Si $Q=x^{2m}+\sum_{j=0}^{2m-1}\alpha_j x^j\in {\mathbb
    F}_2[x]$ est irr\'eductible de degr\'e pair $l=2m$, alors les
  polyn\^omes $Q(x+x^2),Q(1+x+x^2)\in {\mathbb F}_2[x]$ sont
irr\'eductibles de degr\'e $2l$ si la trace
$\hbox{tr}_l(Q)=\alpha_{2m-1}$ de $Q$ vaut $1$ et les deux polyn\^omes
 $Q(x+x^2),Q(1+x+x^2)\in {\mathbb F}_2[x]$ se d\'ecomposent en
un produit de deux polyn\^omes irr\'eductibles de degr\'e $l$ et 
de m\^eme trace sinon.

Si  $Q=x^{2m+1}+\sum_{j=0}^{2m}\alpha_j x^j\in {\mathbb
F}_2[x]$ est irr\'eductible de degr\'e impair $l=2m+1$, alors le
polyn\^ome $Q(1+\alpha_{2m}+x+x^2)\in {\mathbb F}_2[x]$
est irr\'eductible de degr\'e $2l$ et le 
polyn\^ome $Q(\alpha_{2m}+x+x^2)\in {\mathbb F}_2[x]$
se d\'ecompose en un produit de deux polyn\^omes irr\'eductibles
de degr\'e $l$ et de traces diff\'erentes.
\end{thm}

Le th\'eor\`eme \ref{decomp2} est reli\'e \`a une jolie factorisation  
de $x^{2^{2^n}-1}+1\in {\mathbb F}_2[x]$ qui semble nouvelle.
Nous l'appellerons la {\it factorisation de Pascal} car elle fait 
intervenir la r\'eduction modulo $2$ des coefficients binomiaux
${m\choose k}=\frac{m!}{k!(m-k)!}$ constituant le triangle de Pascal.

\begin{thm} \label{Pascalfact} On a 
$$\prod_{k=0}^{2^n-1}\left(1+\sum_{j=0}^{k}{k\choose j}x^{2^j}\right)
=x^{2^{2^n}-1}+1\in {\mathbb F}_2[x]\ .$$
\end{thm}

L'outil principal pour prouver les th\'eor\`emes \ref{decomp2}
et \ref{Pascalfact} est un mono\"{\i}de d\'ecrit dans le chapitre
\ref{monoide}. On l'obtient en consid\'erant
un certain sous-espace vectoriel $\subset{\mathbb F}_p[x]$ qui est
stable pour la composition des polyn\^omes.

\section{Preuves du th\'eor\`eme \ref{tracefact}
et de la proposition \ref{decomp}}

Preuve du th\'eor\`eme \ref{tracefact}: En posant 
$X=\sum_{k=0}^{l-1}x^{p^k}$
dans l'identit\'e triviale $X\prod_{\alpha=1}^{p-1}(X-\alpha)=
X^p-X\in{\mathbb F}_p[X]$ on a 
$$\prod_{\alpha=0}^{p-1}\left(-\alpha+\sum_{k=0}^{l-1}x^{p^k}\right)=
\sum_{k=1}^l x^{p^k}-\sum_{k=0}^{l-1}x^{p^k}=x^{p^l}-x\ .$$
Le polyn\^ome $-\alpha+\sum_{k=0}^{l-1} x^{p^j}$ divise donc
$x^{p^l}-x$ dans ${\mathbb F}_p[x]$.

Par d\'efinition de la trace relative
$\hbox{tr}_l(F)=\sum_{k=0}^{l-1}\rho^{p^k}$ d'un diviseur
irr\'eductible $F$ de $-\alpha+\sum_{k=0}^{l-1}x^{p^k}$, la
racine $\rho\in{\mathbb F}_{p^l}$ de $F$ est \'egalement une racine
de $-\hbox{tr}_l(F)+\sum_{k=0}^{l-1}x^{p^k}$.
Le polyn\^ome irr\'eductible $F\in {\mathbb F}_p[x]$ divise donc   
$-\hbox{tr}_l(F)+\sum_{k=0}^{l-1}x^{p^k}\in {\mathbb F}_p[x]$.
Ceci d\'emontre l'assertion (i).

Posons $q=p^d$ et consid\'erons le corps ${\mathbb F}_q$ \`a $q=p^d$
\'el\'ements. Les arguments utilis\'es dans la preuve de l'assertion
(i) montrent qu'on a \'egalement
$$\prod_{\alpha\in{\mathbb
    F}_q}\left(-\alpha+\sum_{k=0}^{f-1}x^{q^k}\right)=x^{q^f}-x=x^{p^l}-x\in
{\mathbb F}_p[x]\ .$$
La trace $\hbox{tr}_{f,{\mathbb
  F}_q}(F)=\sum_{k=0}^{f-1}\rho^{q^k}\in{\mathbb F}_q$ (o\`u 
$F(\rho)=0$ pour $\rho\in{\mathbb F}_{q^f}$) d'un diviseur
irr\'eductible $F=\sum_j\alpha_jx^j\in{\mathbb F}_q[x]$
de $-\alpha+\sum_{k=0}^{f-1}x^{q^k}$ est donc \'egalement donn\'ee
par $\alpha\in{\mathbb F}_q$. Le produit $\prod_{s=0}^{d-1}
\left(\sum_j \alpha_j^{p^s}x^j\right)$ est une puissance (d'exposant
un diviseur $d$) d'un polyn\^ome irr\'eductible $G\in{\mathbb F}_p[x]$
divisible par $F$. Sa trace
$\hbox{tr}_l(G)=\sum_{k=0}^{l-1}\rho^{p^k}=\sum_{k=0}^{f-1}\left(\rho^{p^{dk}}
+\rho^{p^{dk+1}}+\dots+\rho^{p^{dk+(d-1)}}\right)=\alpha+\alpha^p+\dots+
\alpha^{p^{d-1}}$ vaut donc $d\alpha$ pour $\alpha\in {\mathbb F}_p$.
L'assertion (ii) d\'ecoule maintenant de l'assertion (i).
\hfill$\Box$

Preuve de la proposition \ref{decomp}: 
Soit $\rho$ une racine d'un diviseur
irr\'eductible $F$ de $Q\circ R\in K[x]$. Le corps $K[\rho]=K[x]/(F)$ 
contient
donc le sous-corps $K[R(\rho)]$ qui est une extension de degr\'e $
\hbox{deg}(Q)$
de $K$ car $R(\rho)$ est une racine du polyn\^ome irr\'eductible
$Q\in K[x]$. \hfill $\Box$

\section{Le mono\"{\i}de de composition}\label{monoide}

\begin{lem} \label{lemmonoide}
Pour $A=\epsilon_A+\sum_{k=0}^a\alpha_kx^{p^k}$ et 
$B=\epsilon_B+\sum_{k=0}^b\beta_kx^{p^k}\in {\mathbb F}_p[x]$
avec
$\epsilon_A,\alpha_0,\dots,\alpha_a,\epsilon_B,\beta_0,\dots,\beta_b
\in{\mathbb F}_p$ on a
$$A\circ B=\epsilon_C+\sum_{k=0}^{a+b}\gamma_kx^{p^k}\in{\mathbb
  F}_p[x]$$
avec $\epsilon_C=\epsilon_A+\epsilon_B\sum_{k=0}^a\alpha_k$ et 
$\sum_{k=0}^{a+b}\gamma_kx^k=
\left(\sum_{k=0}^{a}\alpha_kx^k\right)\left(\sum_{k=0}^{b}\beta_kx^k\right)
$.\end{lem}

Preuve: Un calcul facile montre le r\'esultat 
pour le cas particulier $A=\epsilon_A+x^{p^a}$. Le cas g\'en\'eral
s'en d\'eduit par lin\'earit\'e.
\hfill$\Box$

Associons au polyn\^ome 
$A=\epsilon_A+\sum_{k=0}^{a}\alpha_kx^{p^k}\in{\mathbb F}_p[x]$ (avec
$\epsilon_A,\alpha_0,\dots,\alpha_a\in{\mathbb F}_p$) le
symbole $(\epsilon_A,\sum_{k=0}^a \alpha_kx^k)\in {\mathbb
  F}_p\times {\mathbb F}_p[x]$. Par le lemme 
\ref{lemmonoide}, le polyn\^ome compos\'e $A\circ B$ 
de deux tels polyn\^omes de symboles $(\epsilon_A,\alpha)$ et 
$(\epsilon_B,\beta)$ correspond au symbole
$(\epsilon_A+\alpha(1)\epsilon_B,\alpha\beta)$. Ce produit munit 
donc l'ensemble ${\mathcal M}_p={\mathbb
  F}_p\times {\mathbb F}_p[x]$ de ces symboles
d'une structure de mono\"{\i}de associative et distributive \`a gauche
pour la structure d'espace vectoriel \'evidente 
$(\epsilon_A,\alpha)+\lambda(\epsilon_B,\beta)=
(\epsilon_A+\lambda\epsilon_B,\alpha+\lambda\beta)$ sur ${\mathcal M}_p$.
Nous appellerons ${\mathcal M}_p$ le 
{\it mono\"{\i}de de composition}. La projection sur le deuxi\`eme
facteur ${\mathbb F}_p[x]$ de ${\mathcal M}_p$ est un morphisme
de mono\"{\i}de sur le mono\"{\i}de commutatif multiplicatif
${\mathbb F}[x]$. La section $\alpha\longmapsto (0,\alpha)\in{\mathcal
  M}_p$ r\'ealise ${\mathbb F}_p[x]$ comme sous-anneau de l'espace 
vectoriel ${\mathcal M}_p$.
Un autre sous-anneau commutatif (qui n'est cependant pas de type fini)
est d\'efini par le sous-ensemble 
${\mathbb F}_p\times (x-1){\mathbb F}_p[x]$ du mono\"{\i}de
${\mathcal M}_p$, augment\'e de sa structure d'espace vectoriel. 
L'\'el\'ement $(0,1)\in{\mathcal M}_p$ est une
identit\'e bilat\`ere tandis que la multiplication (\`a gauche ou 
\`a droite) par $(0,x)$ correspond \`a l'action de 
l'automorphisme de Frobenius.
Le sous-ensemble ${\mathbb F}_p\times {\mathbb F}_p^*\subset {\mathcal
  M}_p$ est un sous-mono\"{\i}de isomorphe au groupe affine
du corps ${\mathbb F}_p$. Mentionnons \'egalement que le
mono\"{\i}de ${\mathcal M}_p$ poss\`ede des quotients commutatifs
finis de la
forme ${\mathbb F}_p\times ({\mathbb F}_p[x]/(G))$ pour $G\in {\mathbb
  F}_p[x]$ un polyn\^ome divisible par $x-1$. 

\begin{rem} Toutes les d\'efinitions et propri\'et\'es 
\'enonc\'ees dans ce
chapitre restent valable en rempla\c cant $p$ partout par une 
m\^eme puissance
$q=p^e$ d'un nombre premier et en travaillant sur le corps fini
${\mathbb F}_q$ \`a $q$ \'el\'ements.
\end{rem} 

\section{Preuves
des th\'eor\`emes \ref{decomp2} et \ref{Pascalfact}}

Preuve du th\'eor\`eme \ref{decomp2}: 
Consid\'erons d'abord $Q\in{\mathbb F}_2[x]$
irr\'eductible de degr\'e pair $l=2m$.
Comme $(Q_1Q_2)\circ R=(Q_1\circ R)(Q_2\circ R)$, par l'assertion (i)
du th\'eor\`eme \ref{tracefact}   il suffit de montrer
que le polyn\^ome 
$\left(\sum_{k=0}^{2m-1}x^{2^k}\right)\circ(\epsilon+x+x^2)$
divise $x^{2^{2m}}+x$ et que le polyn\^ome
$\left(1+\sum_{k=0}^{2m-1}x^{2^k}\right)\circ(\epsilon+x+x^2)$
est premier \`a $x^{2^{2m}}+x$ (il divisera alors $x^{2^{4m}}+x$
par la proposition \ref{decomp}) pour $\epsilon\in {\mathbb F}_2$. 
En travaillant dans le mono\"{\i}de 
${\mathcal M}_2$, on obtient $(0,\sum_{k=0}^{2m-1}x^k)(\epsilon,1+x)=(
2m\epsilon,1+x^{2m})=(0,1+x^{2m})$ correspondant \`a $x+x^{2^{2m}}$
dans le premier cas. Le deuxi\`eme cas, $(1,\sum_{k=0}^{2m-1}x^k)
(\epsilon,1+x)=(1,1+x^{2m})$, correspond \`a $1+x+x^{2^{2m}}$ qui est
premier \`a $x+x^{2^{2m}}$. La proposition \ref{decomp} 
(ou la factorisation
$(1+x+x^{2^{2m}})(1+(1+x+x^{2^{2m}})^{2^{2m}-1})=
(1+x+x^{2^{2m}})+(1+x+x^{2^{2m}})^{2^{2m}}=x+x^{2^{4m}}$) termine la
preuve pour $Q$ irr\'eductible 
de degr\'e pair $l=2m$. Nous laissons au lecteur le cas similaire
o\`u $Q$ est irr\'eductible de degr\'e impair. \hfill $\Box$

Preuve du th\'eor\`eme \ref{Pascalfact}: Consid\'erons le
polyn\^ome $P_h=1+\sum_{k=0}^h
{h\choose k}x^{2^k}\in{\mathbb F}_2[x]$ correspondant au symbole
$(1,(1+x)^h)\in{\mathcal M}_2$. Les identit\'es faciles
$(1,1+x)(1,(1+x)^h)=(1,(1+x)^h)(0,1+x)=(1,(1+x)^{h+1})\in{\mathcal
  M}_2$ pour $h\in{\mathbb N}$ montrent 
qu'on a $P_1\circ P_h=P_h\circ(x+x^2)=P_{h+1}$ pour
tout $h\in{\mathbb N}$. En it\'erant l'identit\'e
$P_1\circ P_h=1+P_h+P_h^2=P_{h+1}$ on obtient
$$\begin{array}{lcl}
P_{h+1}&=&1+P_h(1+P_h)=1+P_hP_{h-1}(1+P_{h-1})=\dots\\
&=&1+P_hP_{h-1}\dots P_1P_0(1+P_0)=1+x\prod_{k=0}^h P_k\ .\end{array}$$
Pour $h+1=2^n$, on a donc l'\'egalit\'e
$$P_{2^n}=1+x+x^{2^{2^n}}=1+x\prod_{k=0}^{2^n-1} P_k$$
qui d\'emontre le th\'eor\`eme.\hfill $\Box$

\begin{rem} Le produit 
$$P_{2^{n-1}}\cdots
  P_{2^n-1}=\frac{x^{2^{2^n}-1}+1}{x^{2^{2^{n-1}}-1}+1}\in
{\mathbb F}_2[x]$$
s'identifie au produit de tous les polyn\^omes irr\'eductibles 
distincts de degr\'e
$2^n$ dans ${\mathbb F}_2[x]$. Le th\'eor\`eme \ref{tracefact}
appliqu\'e \`a $P_{2^{n-1}}=1+\sum_{k=0}^{{2^n}-1}x^{2^k}$ (ou le
th\'eor\`eme \ref{decomp2} appliqu\'e aux identit\'es 
$P_{h+1}=P_h\circ(x+x^2)$) montre qu'un tel polyn\^ome irr\'eductible 
$F=\sum_{k=0}^{2^n}\alpha_kx^k$ divise le dernier facteur
$P_{2^{n-1}}$ si et seulement si sa trace
$\hbox{tr}_{2^n}(F)=\alpha_{2^n-1}$ vaut $1$. Mentionnons \'egalement
que les formules 
$P_{h+1}=P_h\circ P_1=P_h\circ(x+x^2)$ pour $h\geq 1$ illustrent et
pr\'ecisent le th\'eor\`eme \ref{decomp2}. 
\end{rem}

\begin{rem} Pour $p$ un premier quelconque (ou plus g\'en\'eralement 
pour $q=p^e$ une puissance d'un premier), les polyn\^omes
$$P_{h,\alpha}=-\alpha+\sum_{k=0}^h{h\choose k}x^{p^k}\in{\mathbb F}_p[x]$$ 
ont \'egalement des propri\'et\'es int\'eressantes: 
$P_{h,\alpha}$ divise $P_{h+1,2\alpha}$ car
$$\begin{array}{lcl}0&=&\left(-\alpha+\sum_k{h\choose k}\rho^{p^k}\right)+
\left(-\alpha+\sum_k{h\choose k}\rho^{p^k}\right)^p\\
&=&
-\alpha+\sum_k{h\choose k}\rho^{p^k}
-\alpha+\sum_k{h\choose k-1}\rho^{p^k}=P_{h+1,2\alpha}\end{array}$$
pour $\rho\in\overline{\mathbb F}_p$ une racine de $P_{h,\alpha}$.
Le polyn\^ome $P_{h,\alpha}$ divise donc $P_{p^n,2^{p^n-h}\alpha}$
pour $h\leq p^n$. Il divise donc \'egalement $x^{p^{2p^n}}-x$
en appliquant le th\'eor\`eme \ref{tracefact} \`a
$P_{p^n,2^{p^n-h}\alpha}=2^{p^n-h}\alpha+x+x^{p^{p^n}}$.
L'ensemble des polyn\^omes $P_{h,\alpha}$ correspond au 
sous-mono\"{\i}de ${\mathbb F}_p\times (1+x)^{\mathbb N}$ de 
${\mathcal M}_p$. Il est donc ferm\'e pour la composition
et on a $P_{h,\alpha}\circ P_{m,\beta}=P_{h+m,\alpha+2^h\beta}$.
\end{rem}

Institut Fourier,
Laboratoire de Math\'ematiques,
UMR 5582 (UJF-CNRS),
100, rue des Math\'ematiques,
BP 74,
38402 St MARTIN D'H\`ERES Cedex, France

Adresse courriel: Roland.Bacher@ujf-grenoble.fr
\end{document}